\magnification\magstep1
 \input epsf.sty
\input  amssym.tex
\input color

 \centerline{\bf Invertible Dynamics on Blow-ups of ${\Bbb P}^k$}

 \medskip
 
 \centerline{Eric Bedford}
 
 \medskip
 \medskip\noindent{\bf \S0.  Introduction }  Let $X$ be a complex manifold, and let $f$ be an automorphism, i.e. biholomorphic self map, of $X$.  We discuss the project of finding compact complex manifolds $X$ which carry automorphisms $f$ which are dynamically interesting.  If the dimension of $X$ is 1, then $X$ is a compact Riemann surface, and if the genus is $\ge2$, then $Aut(X)$ is finite.  The other two cases are when $X$ is a torus, in which case the automorphisms are essentially translations, or $X={\Bbb P}^1$, in which case the automorphisms are linear (fractional).  Thus, to find maps with interesting dynamics, we must start with dimension 2.  A Theorem of Cantat (see [C1,2]) restricts the set of possible surfaces $X$.  However, little is known about exactly which surfaces $X$ might arise in these cases.  Some basic ergodic properties hold for all automorphisms with dynamical degree $>1$, but little is known about the topological properties of such maps.  Here we focus on the case where $X$ is a blowup of ${\Bbb P}^2$ and we describe some of the results that are known.  
  
Next we will consider the case of dimension 3 and higher.  In this case, we find that it is natural to widen our search to manifolds $X$ which carry pseudo-automorphisms.  These are birational maps whose indeterminate and exceptional behaviors only influence subvarieties of codimension 2 and greater.  
It is the purpose of this paper to discuss pseudo-automorphisms that can be obtained from birational maps of ${\Bbb P}^k$ by blowing up, and in this discussion we formulate a number of open questions.

\bigskip\noindent {\bf \S1. Rational surfaces }    The group of automorphisms (biholomorphic self-maps) of complex projective space ${\Bbb P}^k$ is $PGL({\Bbb C},k+1)$, which may be written as linear fractional transformations of ${\Bbb C}^k$.  We can widen this class of manifolds  by blowing up.  In terms of complex analysis, all the global holomorphic functions are constant, and this is not changed by blowing up.  However, the set of biholomorphic self-mappings might change.

Let us start with dimension 2 and the classical example of the Cremona Involution, which on ${\Bbb C}^2$ is given by $(x,y)\mapsto (1/x,1/y)$, and on projective space we write it as a mapping of degree 2:
$$J: [x_0:x_1:x_2] \mapsto [1/x_0:1/x_1:1/x_2] = [x_1x_2:x_0x_2:x_0x_1]$$

A rational (or meromorphic) mapping is said to be {\it regular} at a point $p$ if it is holomorphic in a  neighborhood of $p$.  The {\it indeterminacy locus} of a rational map $f:X\dasharrow Y$, written ${\cal I}(f)$ is defined as the set of all points where $f$ is not regular.  It may be shown that if $f$ is indeterminate at a point $p$, then $f$ blows up $p$ to a variety $V\subset Y$, and the dimension of $V$ is at least one.  This ``blow up'' image $V$ may be defined in more than one way. One of them is essentially the cluster set:
$V=\bigcap_{\epsilon>0}  {\rm closure}\left (f(B(p,\epsilon)-{\cal I}(f)) \right)$,
where  $B(p,\epsilon)$ denotes the ball about $p$ with radius $\epsilon$.

If $W\subset X$ is any subvariety, and if ${\cal I}(f)$ does not contain any irreducible component of $W$, then $W-{\cal I}(f)$ is dense in $W$.  The closure of the image $f(W-{\cal I}(f))$ is a subvariety of $Y$, and we call it the {\it strict transform} of $W$.
A general fact is that the indeterminacy locus has codimension $\ge2$.  Thus if $H\subset X$ then we may take its strict transform $f(H)\subset Y$.  We say that $H$ is {\it exceptional} if the codimension of $f(H)$ is $\ge2$.

With these definitions, we see that the indeterminacy locus of $J$ is ${\cal I}(J)=\{e_0,e_1,e_2\}$, and the lines $\Sigma_j:=\{x_j=0\}$, $j=0,1,2$, are exceptional.  Since $J$ is an involution, (i.e.\ $J^2=$identity), we see that $J$ blows up $e_j$ to $\Sigma_j$.
Geometrically, $J$ acts as an inversion in the coordinate triangle, as shown in Figure 1.

\medskip
\epsfysize=.85in
\centerline{ \epsfbox{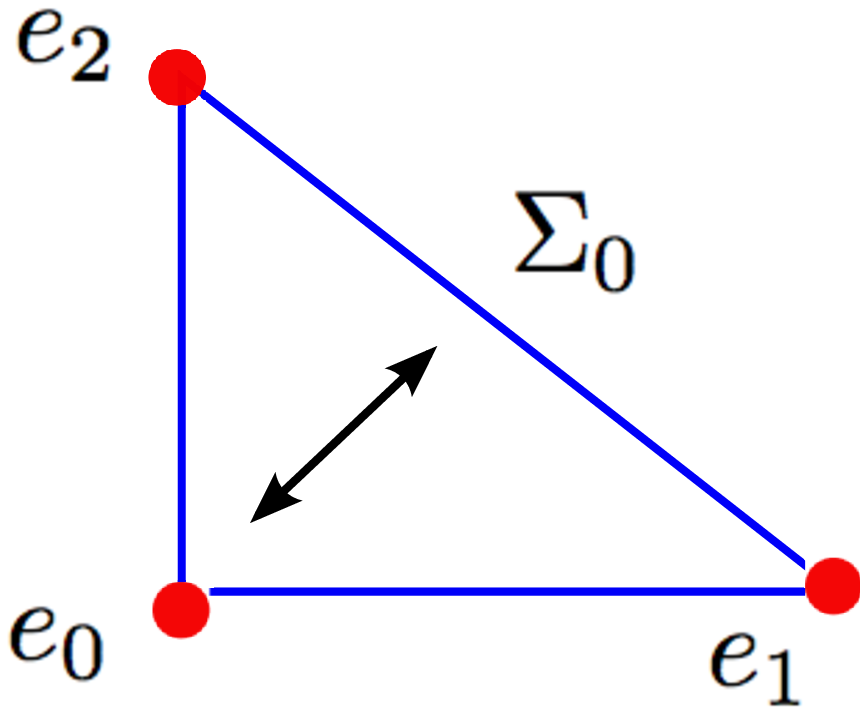} }

\centerline{ indeterminate\ point $e_0=[1:0:0] \leftrightarrow$ exceptional curve $\Sigma_0=\{x_0=0\}$}

\centerline{Figure 1. }

\medskip

Now we blow up of ${\Bbb P}^2$ at the point $e_0$.  This is a new manifold $X$ with a holomorphic projection $\pi:X\to{\Bbb P}^2$ with the properties:  (1)  the exceptional fiber $E:=\pi^{-1}(e_0)$ is isomorphic to ${\Bbb P}^1$, and (2) $\pi:X-E\to {\Bbb P}^2-\{e_0\}$ is biholomorphic.  We may represent $X$ in local coordinates  $(\xi_1,\xi_2)$ over ${\Bbb P}^2-\Sigma_0$ 
$$ \ \pi(\xi_1,\xi_2) = [1: \xi_1:\xi_1\xi_2] = [x_0:x_1:x_2],  \ \pi^{-1}(x) = (\xi_1=x_1/x_0, \xi_2= x_2/x_1)$$
\medskip
\epsfysize=1.35in
\centerline{\epsfbox{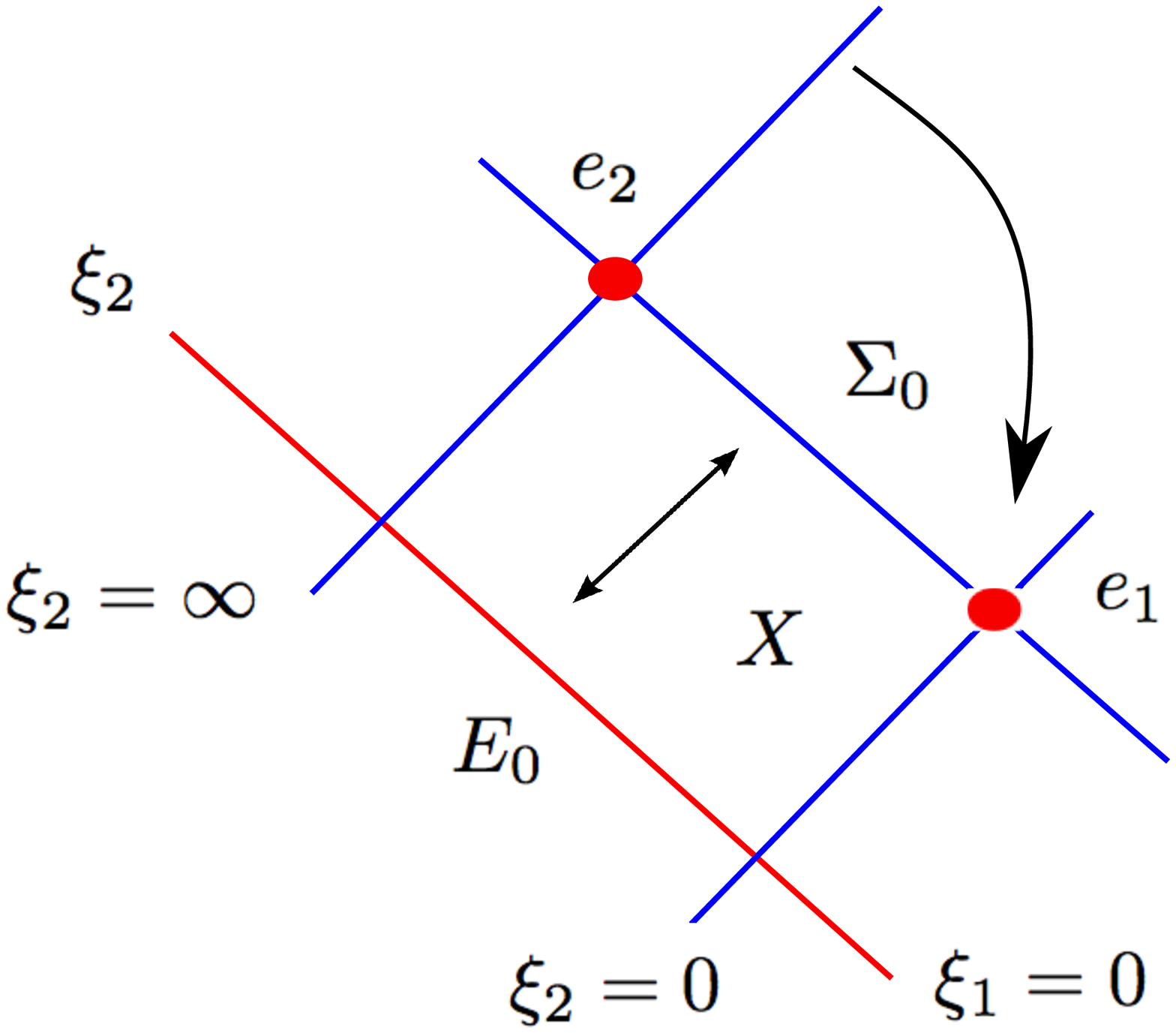}}

\centerline{Figure 2.}
\smallskip
Figure 2 shows the new blowup space $X$ in the $(\xi_1,\xi_2)$ coordinate chart.  The point $e_0$ has been replaced by a curve (called {\it exceptional} or {\it blowup divisor}) $E_0$.   The projection $\pi$ maps  $X-E_0$ biholomorphically to ${\Bbb P}^2-\{e_0\}$.  Since $e_0\notin \Sigma_0$, $\pi^{-1}$ is holomorphic in a neighborhood of $\Sigma_0$, and via this biholomorphism, we have a new curve $\pi^{-1}\Sigma_0$ inside $X$.  We write this again as $\Sigma_0$, although technically it is called the {\it strict transform} of $\Sigma_0$ in $X$.  We may also use $\pi^{-1}$ to lift the curve $\Sigma_1-\{e_0\}$ to $X-E_0$; inside the $(\xi_1,\xi_2)$ coordinate chart,  this corresponds $\xi_2=\infty$.  The closure of this set is a curve in $X$, which is the strict transform $\Sigma_1'$ of $\Sigma_1$.  Although the curves $\Sigma_1$ and $\Sigma_1'$ are isomorphic, their normal bundles are not.

Now we describe the behavior of the induced map $J_X:=\pi^{-1}\circ J\circ\pi$ in the part of the $(\xi_1,\xi_2)$-coordinate chart  where $\xi_2\ne0,\infty$.  $J_X$ will map this set into $X-E_0$, which is mapped biholomorphically to ${\Bbb P}^2-\{e_0\}$ by $\pi$.  We represent its range in the homogeneous $x$-coordinates as
$$J_X(\xi_1,\xi_2) = J([1:\xi_1:\xi_1\xi_2])= [1:\xi_1^{-1}:\xi_1^{-1}\xi_2^{-1}]=[\xi_1\xi_2:\xi_2:1]$$   
which shows that $J_X$ is regular in this coordinate chart.

To summarize: the induced map $J_X$ has two points of indeterminacy $e_1$ and $e_2$.  The exceptional locus of $J_X$ consists of the strict transforms of $\Sigma_1$ and $\Sigma_2$.  If we continue this process and blowup all 3 points $e_0$, $e_1$ and $e_2$, then we obtain a manifold $Z$, and the induced map $J_Z$ is an automorphism of $Z$.  This is pictured symbolically in Figure 3, where all 3 blowup divisors are represented as red segments, and the strict transforms of the $\Sigma_j$ are blue segments.  Together, they form a hexagon, and the action of the induced map $J_Z$ is to interchange opposite sides of the hexagon.

\medskip
\epsfysize=1.25in
\centerline{\epsfbox{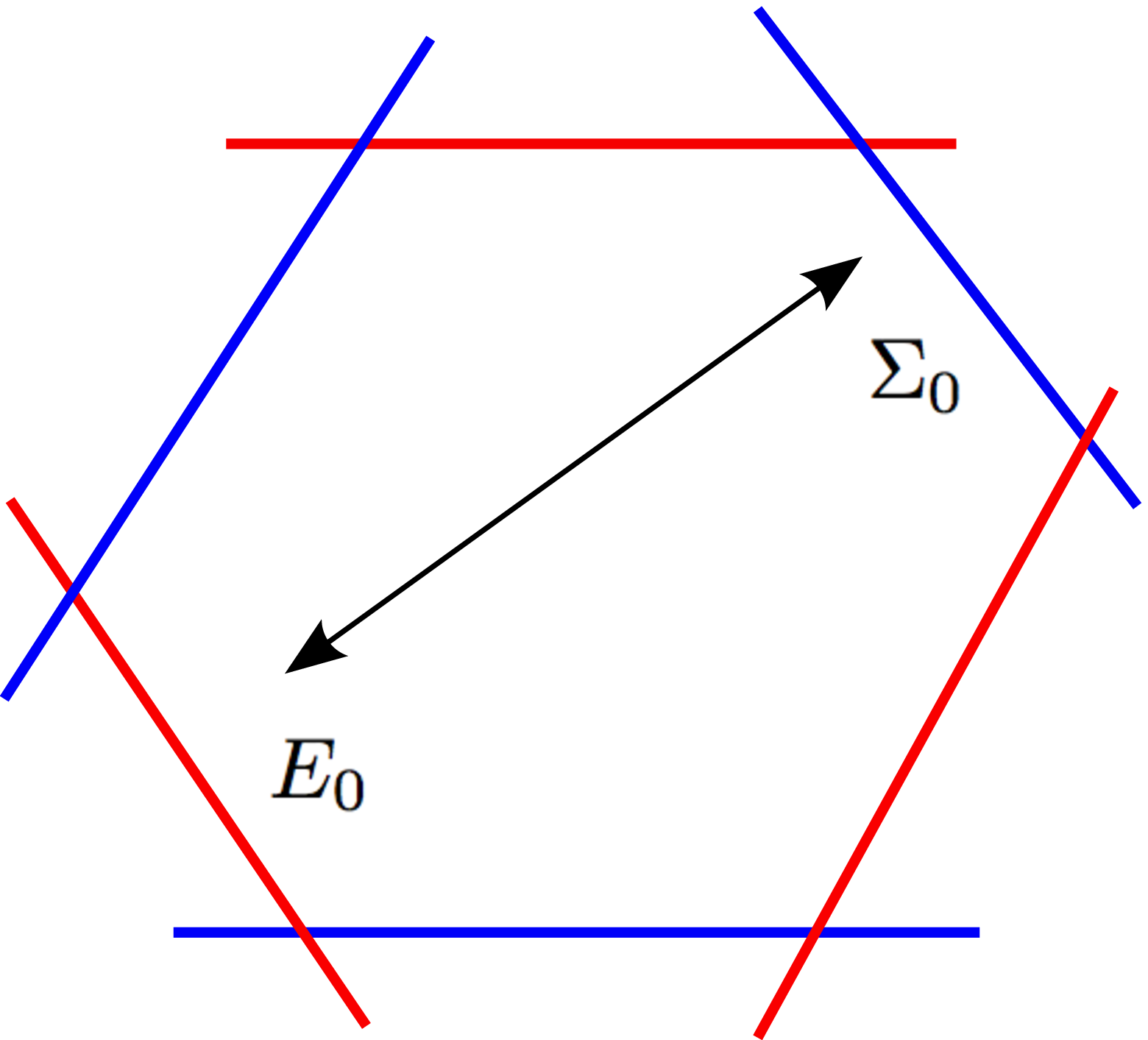} }

\centerline{Figure 3.}
\smallskip

Thus we have started with a birational map and  have obtained an automorphism on some blowup of ${\Bbb P}^2$.  We will explore the question of how much more generally this might work:  what are the birational maps of ${\Bbb P}^2$ that might lead to automorphisms after some blowups?  For instance, de Fernex and Ein [dFE] have shown:  {\sl  If $f:X\dasharrow X$ is a birational map of finite order (i.e.\ if $f^N$ is the identity for some $N$), then there is an iterated blowup $\pi: Z\to X$ such that the induced map $f_Z$ is an automorphism of $Z$.}

\medskip\noindent{\bf \S2.  Degree complexity, or dynamical degree }  Every holomorphic or rational map $f:{\Bbb P}^k\dasharrow {\Bbb P}^k$ is represented by polynomials $f=[f_0:\cdots:f_k]$ of a common degree $d$.  Dividing by the GCD, we may suppose that the degree $d$ is minimal, and we set ${\rm deg}(f):= d$.  The {\it dynamical degree} is the limit
$$\delta(f) := \lim_{n\to\infty} \left( {\rm deg} (f^n) \right)^{1/n}.$$
We may think of $\delta$ as the degree complexity, or as  the exponential rate of degree growth of $f^n$ as $n\to\infty$.  We note that $1\le \delta\le d$, and for ``generic'' $f$ we have $\delta=d$ (see [FS]).   The cases we are interested in, however, are when $f$ is an automorphism (or pseudoautomorphism) with $\delta>1$, and in this case $\delta$ is algebraic but never rational (see [B]).  In particular, it follows that $\delta<d$.  

We may consider the lift $F=(F_0,\dots,F_k)$ of $f$ to a polynomial self-map of ${\Bbb C}^{k+1}$.  In this case, $F^n=F\circ\cdots\circ F$ is the usual composition of polynomials, and the degree of  $F^n$ is $d^n$.  Let us write  $F^n=\psi_n F^{\langle n\rangle}$, where $\psi_n$ is the GCD of $F^n$.  In this case, we have ${\rm deg}(F^{\langle n\rangle})\sim \delta^n$, and ${\rm deg}(\psi_n)\sim d^n-\delta^n$, so $\psi_n$ carries almost all the degree of $F^n$.

The degree of $f$ is closely related to how it pulls back hypersurfaces.  Recall that $H^2({\Bbb P}^k)=H^{1,1}({\Bbb P}^k)$, and $H^2({\Bbb P}^k;{\Bbb Z})$ is generated by the class of a hyperplane, which we will again write as $H$.   if $V$ is any hypersurface of degree $m$, then its class is given by ${\rm deg}(V)\cdot H\in H^2({\Bbb P}^k;{\Bbb Z})$.  The pullback of the class of a hyperplane $H=\{\sum c_jx_j=0\}$ is given by the class of $\{\sum c_jf_j=0\}$.  Thus $f^*H={\rm deg}(f) H$.

Let $\beta$ denote any K\"ahler form (for instance, the Fubini-Study form) with $\int_{{\Bbb P}^k}\beta^k=1$.  Then ${\rm deg}(V)=\int_V\beta^{k-1}$.  It follows that for a generic hypersurface $H$, we have
$$\eqalign{ \delta = &\lim_{n\to\infty}\left( {\rm Vol}_{k-1}(f^{-n}(H))\right)^{1/n}  \cr
&=\lim_{n\to\infty}\left(  \int_{f^{-n}H}\beta^{k-1} \right)^{1/n}  = \lim_{n\to\infty} \left( \int_{{\Bbb P}^k}\beta^{k-1}\wedge (f^n)^*\beta\right)^{1/n} }$$
so the dynamical degree also measures the exponential rate of growth of $(k-1)$-dimensional volume under pullback.

It would be convenient if we could have  $(f^n)^*=(f^*)^n$.  In our case, that would mean that ${\rm deg}(f^n)=({\rm deg}(f))^n$  (see [FS]).   In dimension $k=2$, [DF] showed that there is an iterated blowup $\pi:X\to{\Bbb P}^2$ such that the induced map $f_X:= \pi^{-1}\circ f\circ \pi$ does satisfy    $(f_X^n)^*=(f_X^*)^n$.   In this case, we let $\beta_X:=\pi^*(\beta)$, so these integrals become:
$$ \eqalign{\delta= &\lim_{n\to\infty} \left( \int_{{\Bbb P}^2} \beta\wedge (f^n)^*\beta \right)^{1/n} = 
\lim_{n\to\infty}\left(\int_{X} \beta_X\wedge (f^n_X)^*\beta_X \right)^{1/n} 
\cr  &=\lim_{n\to\infty}\left(\int_{X}\beta_X\wedge (f_X^*)^n\beta_X \right)^{1/n}= ||f^*_X||_{\rm sp}}$$
where $||\cdot||_{\rm sp}$ denote the spectral radius, i.e.\ the modulus of the largest eigenvalue.  The reason that the growth of $\beta_X$ under $(f^*_X)^n$ gives the growth of $||(f_X^n)^*||$ is that since it is a K\"ahler form, it can lie in an eigenspace only if $\delta=1$, and thus $(f^*_X)^n\beta_X$ must grow like the largest eigenvalue.

To give a simple example of the regularization of a map, we return to the map $J$ from \S1.  In this case, we have $J^*=2$, which means that $J^*$ acts on $H^2$ by multiplication by 2.  On the other hand $J^2$ is the identity map, so $(J^2)^*=1\ne 4 = (J^*)^2$.

Now consider the space $X$, which was obtained by blowing up ${\Bbb P}^2$ at the points $e_0$, $e_1$, $e_2$.  The cohomology group $H^2(X;{\Bbb Z})$ has the ordered basis $\langle H_X,E_0,E_1,E_2\rangle$, where $H$ denotes the class of the strict transform of a generic hypersurface (line), and $E_j$ denotes the class of the exceptional fiber over $e_j$.  Let us see how to represent $J_X^*$ with respect to this basis.   Since $J_X$ is an automorphism of $X$, we see that $J_X^*: E_j\leftrightarrow \Sigma_j'$.  Next we need to represent $\Sigma_j'$ with respect to this basis.  Inside ${\Bbb P}^2$, we have that $\Sigma_j=H$ is the generator of $H^2({\Bbb P}^2;{\Bbb Z})$.  Now let us pull this back under $\pi^*$.  It follows that $\pi^*H=H_X=\pi^*\Sigma_j$.  On the other hand, suppose for instance that $j=0$.  Then $\Sigma_0$ contains $e_1$ and $e_2$, so when we pull back, we get the total preimage, and $\pi^*(\Sigma_0)=\Sigma_0' + E_1+E_2$.  This gives $J_X^*: E_j\mapsto \Sigma_j'= H_X-E_1-E_2$.  Finally, $J^{-1}\{\sum c_jx_j=0\} =\{c_0 x_1x_2+c_1 x_0x_2 + c_2 x_0x_1=0\}$, so $J^{-1}\{\sum c_jx_j=0\} $ contains $\{e_0,e_1,e_2\}$. This is a curve of degree 2, so its class in $H^2({\Bbb P}^2)$ is $2H$.   Applying $\pi^*$ once again, we see that inside $H^2(X)$ we have $2H_X = J^{-1}\{\sum c_j x_j=0\}_X + \sum E_p$.  Since $J_X$ is an automorphism, we obtain $J^*_X= J^{-1}\{\sum c_j x_j=0\}_X = 2H - \sum E_p$.  Writing this as a matrix with respect to our ordered basis, we find
$J^*_X = \pmatrix {2 & 1& 1& 1\cr  -1& 0 & -1 &-1\cr -1& -1& 0 & -1\cr  -1& -1 & -1 &0}$, which satisfies $(J_X^*)^2=I$. 

In the preceding heuristic argument we did not take into account the multiplicities of the $E_p$, which in fact turn out to be 1 in this case.  In higher dimension, they are $>1$.

\bigskip\noindent {\bf  \S3.  Some rational surface automorphisms with $\delta>1$.}
We start with a particularly simple family of planar rational maps:  
$$f_{a,b}(x,y)=\left(y,{\displaystyle y+a \over \displaystyle x+b}\right):{\Bbb C}^2\dasharrow{\Bbb C}^2$$  
We may write this in homogeneous coordinates $[x_0:x_1:x_2] = [1:x:y]$ as
$$f_{a,b}: \ \ [x_0:x_1:x_2] \mapsto [x_0(b x_0+ x_1): x_2(b x_0+ x_1): x_0(ax_0 + x_2)]$$
The points of indeterminacy are $\{e_1,e_2,p=[1:-b:-a]\}$.  The Jacobian determinant is $2 x_0 (b x_0 + x_1) (a x_0 + x_2)$, which vanishes on three lines, which are the exceptional curves of $f_{a,b}$.  Choosing a linear map mapping this triangle to the standard coordinate triangle $\{x_0x_1x_2=0\}$, we may conjugate $f_{a,b}$  to a map
of the form $L\circ J$, where $L$ is linear, and $J$ is the map from \S1.
The triangle of exceptional curves, shown in Figure 4, is mapped as:  
$$\Sigma_\beta=\{bx_0+x_1=0\}\to e_2 \dasharrow \Sigma_0\to e_1\dasharrow \Sigma_\gamma=\{ax_0+x_2=0\}\to q:=[1:-a:0]$$

Now let us construct a new space $Y$ by blowing up ${\Bbb P}^2$ at the points $e_1$ and $e_2$.  We see that the induced map $f_Y$ has one exceptional curve $\Sigma_\gamma$ and one point of indeterminacy $p$.  If the orbit of $q$ lands on $p$, i.e.\ if $f^Nq=p$, then we may blow up the orbit $q$, $fq$, \dots, $f^Nq=p$.  If we start by blowing up $q$, then $\Sigma_\gamma$ is no longer exceptional, but the blowup fiber $Q$ over $q$ is exceptional.  Next, we blow up $fq$, so $Q$ is now mapped to $fQ$ (the fiber over $q$) and is no longer exceptional.  We continue over the whole orbit and obtain an automorphism, since the fiber $P$ over the indeterminate point $p$ now maps in a regular way to the line $\overline{q,e_2}$.  (It is easy to see that this line is what we would expect because $p=\Sigma_\gamma\cap \Sigma_\beta$, and $\Sigma_\gamma\mapsto q$,  $\Sigma_\beta\mapsto e_2$.)

\bigskip
\epsfxsize=2.25in
\centerline{\epsfbox  {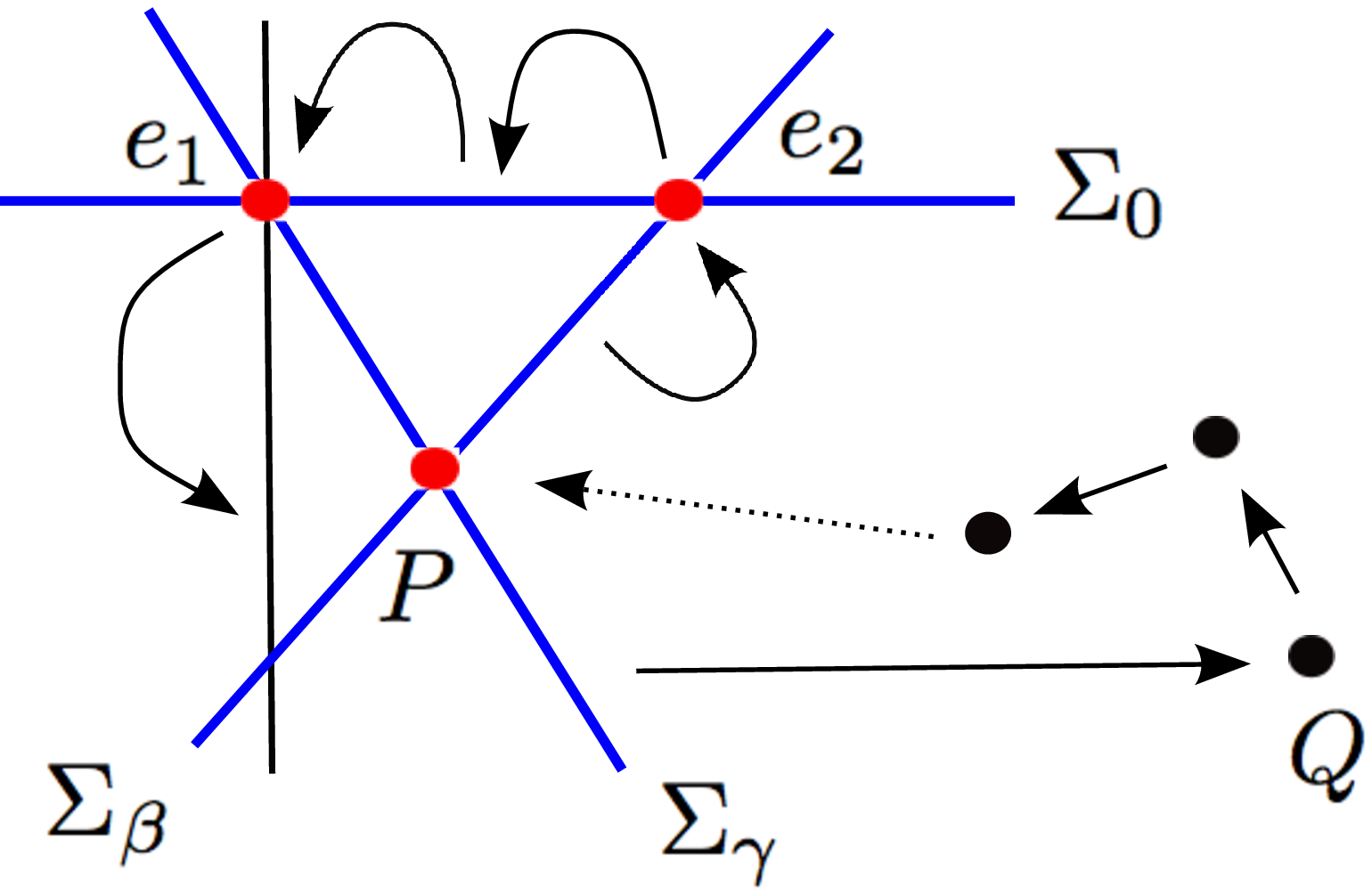}  }

\centerline{Figure 4.}
\bigskip

Let us define the set 
$${\cal V}_N:=\{(a,b)\in{\Bbb C}^2: f^N_{a,b}(q) = f^N_{a,b}(-a,0) = p = (-b,-a)\}$$  
From [BK2, Theorem 2] we have:  {\it  If $(a,b)\in{\cal V}_N$, then there is a blowup $X$ such that the induced map $f_X$ is an automorphism.}  
 
 Now let us suppose that $f_{a,b}$ is an automorphism and show how determine $f^*_X$.  We may use the basis $\langle H_X,E_1,E_2,P=f^NQ, f^{N-1}Q,\dots,Q\rangle$ as an ordered basis for $H^2(X;{\Bbb Z})$.  Since $f_X$ is an automorphism, we can read off the behavior of $f^{-1}_X$ from Figure 4:
 $$f^*_X: \ \  E_1\mapsto \Sigma_0',\  E_2\mapsto \Sigma_\beta',\  P=f^NQ\mapsto f^{N-1}Q\mapsto\cdots\mapsto Q\mapsto \Sigma_\gamma'$$
Further, as we saw in the case of $J_X$, we have
$$\Sigma_0'=H_X - E_1-E_2, \ \ \Sigma_\beta'=H_X-P-E_2,\  \ \Sigma_\gamma'=H_X-E_1-P$$
and
$$f^*_X: \ \ H_X \mapsto 2H_X - E_1-E_2-P$$
This defines $f_X^*$ on all the basis elements of $H^2(X)$, so we may now compute the characteristic polynomial of $f^*_X$, and we find that it is:  
$$\chi_N(t) = t^{N+1}(t^3-t-1) + t^3+t^2-1 \eqno(3.1)$$
The dynamical degree of $f_{a,b}$ is the spectral radius of $f^*_X$, which is the largest root of $\chi_N$.  We conclude  (cf [BK2]) that: {\it If $(a,b)\in {\cal V}_N$, and if $N\ge 7$, then $f_{a,b}$ is an automorphism of $X$ with $\delta(f_{a,b})>1$.} 

It is known that ${\cal V}_N\ne\emptyset$ for all $N$ (see [M1] and [BK3]), but: {\it It is not known whether ${\cal V}_N$ is discrete for $N\ge7$.}  For generic $(a,b)\in{\Bbb C}^2$, the dynamical degree of $f_{a,b}$ is $\delta_*\sim 1.324$, the largest root of $t^3 - t  - 1$. {\it Does the cardinality of ${\cal V}_N$ grow like $\delta_*^N$ as $N\to\infty$?}

\bigskip\noindent{\bf \S4.  Connection between dynamical degree and length growth: A graphic example }  
If we choose $(a,b)\in{\cal V}_7$, then the map $f_{a,b}$ in \S3 will be an automorphism with $\delta\sim1.17628$, which is the largest root of $\chi_7$ from (3.1).  By the discussion in \S2, we know that the $n$th iterate of a (complex) line will have 2-dimensional area $\sim\delta^n$.    This is closely related to the fact (see [C]) that the entropy of $f_{a,b}$ is $\log\delta$.  There is only one map $f_{a,b}$ (and its inverse) with $(a,b)\in{\cal V}_7\cap {\Bbb R}^2$.  The restriction of this map to $X_R$, the real points of $X$ was shown in [BK3] to have entropy $\log\delta$.   The picture below, taken from [BK2, Figure~A.1], shows an example of a real point $(a,b)\in{\cal V}_7$.   This shows the image of a line $L$ after $n$ iterations, and inspection shows empirically that the length grows $\sim\delta^n$.

To represent the real projective plane, we have taken the usual polar coordinates $(r,\theta)$ and replaced them by modified polar coordinates $(\rho,\theta)$, where $\rho=\arctan(r)$.  Thus a linear hypersurface will not appear straight in these coordinates.  The green curve is an invariant cubic, and the red points (on the invariant curve) are the two fixed points.  The three exceptional curves $\Sigma_0$, $\Sigma_\beta$ and $\Sigma_\gamma$ are labeled.  The intersection points of any two of these exceptional curves are necessarily indeterminate.  The image of $\Sigma_\gamma$ is labeled ``0'', and the image of ``0'' is ``1'', etc.  The point ``7'' is the indeterminate point $\Sigma_\beta\cap\Sigma_\gamma$.

The forward image of the line $L$ appears to be ``bunched'' at the points ``0'', ``1'', \dots\ \ These points are blown up in the construction of $X$, so the ``bunching'' is an artifact of the projection $\pi$, which takes all the points of a blowup fiber and collapses them to a point.  It appears that $f^nL$ may converge to a lamination as $n\to\infty$, and it would be interesting to know whether this is true.  (A related map was shown to have a lamination in [BD].)

\medskip
\epsfxsize=2.65in
\centerline{\epsfbox{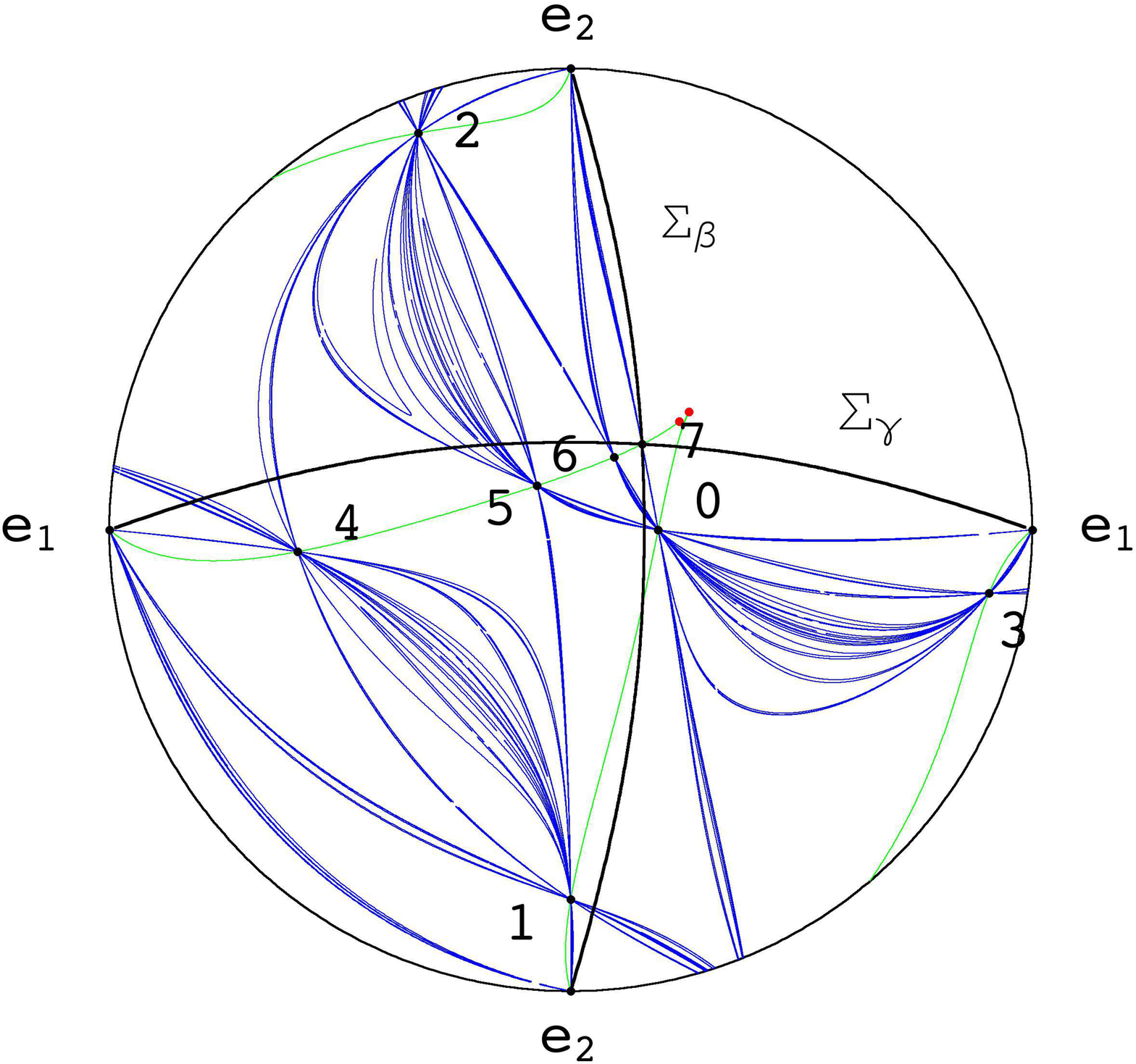}  }

\centerline{$f(x,y) = \left(y,{\displaystyle y+a \over \displaystyle x+b} \right)$, $a= - 0.499497$, $b=-0.415761$}

\centerline{Figure 5.}

We note that one of the fixed points is attracting, and the basin has full area  (see [M1], [BK3]).  Another  graphical representation for this map is to choose a point of the basin and plot its orbit under $f^{-1}$, as was done in [M1, Figure 1].  It is striking to see the similarities between these figures, the main visual difference coming from the coordinate systems  used:  affine for [M1, Figure~1] and  polar for Figure 5 above.

\bigskip\noindent{\bf \S5.  Heuristic picture:  Dynamical complexity vs 
degree complexity}  The main goal of this paper is to discuss maps with interesting dynamics, so let us give some examples of interesting dynamics.  In dimension 1, we consider a rational map $f:{\Bbb P}^1\to {\Bbb P}^1$ with 
degree$(f) =\delta(f) = {\rm deg}_{\rm top}(f)>1$.  Such maps are not invertible.
One of the basic results of the subject concerns the backward dynamics of $f$, i.e. the distribution of preimages of a point.

\medskip
\epsfxsize=2.05in
\centerline{\epsfbox{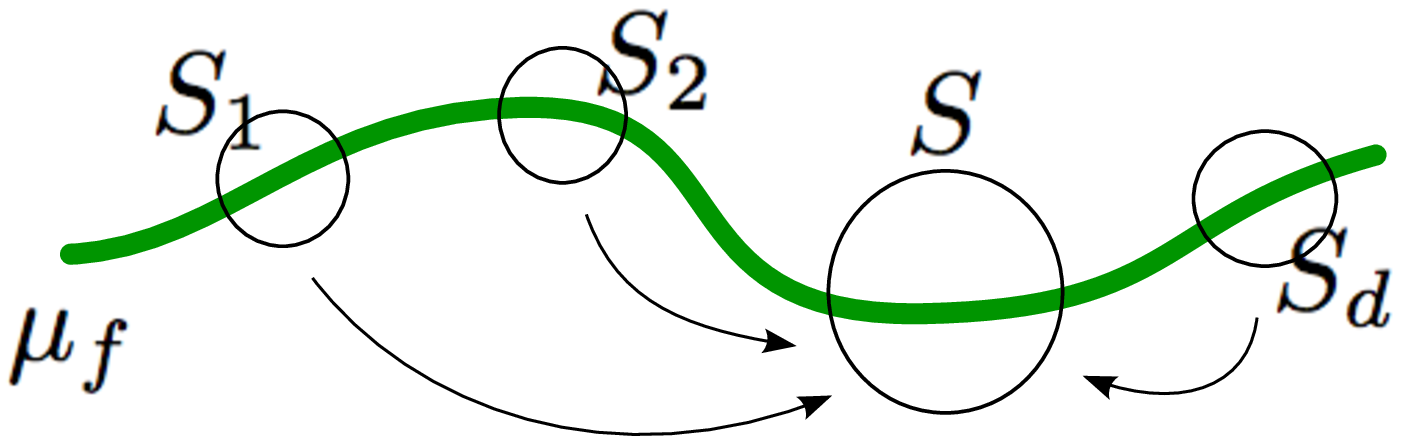} }

\centerline{Choosing preimages of a point is like flipping a $d$-sided coin.}

\centerline{Figure 6.}

\proclaim Theorem ([Br], [FLM], [L]).  For almost all  points $z_0$ there is a limiting distribution of point masses over the preimages of $z_0$: 
$$\mu_f:=\lim_{n\to\infty} {1 \over d^n} \sum_{\{a:f^n(a) = z_0\} }\delta_a$$

This measure is {\it balanced}, which means that, locally, $f^*\mu_f = d\cdot\mu_f$.  Thus, as we consider backward iteration,  the different branches of $f^{-1}$ give  $\mu_f(S_j) = \mu_f(S)/d$.   The effect is like Bernoulli trials, as illustrated in Figure 6.

If $K\subset X$ is a compact set, and if $f:X\to X$, then we may define the stable set $W^s(K):=\{x\in X:{\rm dist}(f^n(x),K)\to 0{\rm\ as\ }n\to\infty\}$.
In dimension 1, the stable set of a point $W^s(x_0)$ is just the set of all preimages of $x_0$.  (The unstable set $W^u(x_0)$ is defined as above, with $f$ replaced by $f^{-1}$.) The Theorem above says that the asymptotic distribution of $W^s(x_0)$ is independent of $x_0$.  In dimension 2, the stable set has the structure of a manifold.
\proclaim Stable Manifold Theorem. Let  $f:X\to X$ be an automorphism of a complex surface, let $x_0$ be a saddle fixed point, and let
$$W^s(x_0):=\{x: \lim_{n\to\infty}{\rm dist}(f^nx,x_0)=0\}$$
be its stable set.  Then there is an injective holomorphic immersion $\xi:{\Bbb C}\to X$ such that $\xi({\Bbb C})=W^s(x_0)$.

A classic example is the Horseshoe Map.   Figure 7  is the ${\Bbb R}^2$ slice of a complex automorphism of ${\Bbb C}^2$.  It shows the stable and unstable manifolds of a saddle fixed point $x_0$, which is marked in the upper left.  The arcs which are essentially ``left-right'' (with bends) are contained in $W^u(x_0)$.  They are cut off by the viewbox of the picture, but they are connected in ${\Bbb R}^2$.  Similarly, the arcs which are ``up-down'' oriented are contained in the stable manifold $W^s(x_0)$.  In the case of the horseshoe, the closure of $W^s(x_0)$ is a complicated set; a laminar structure is clearly visible.

\medskip
\epsfxsize=2.05in
\centerline{\epsfbox{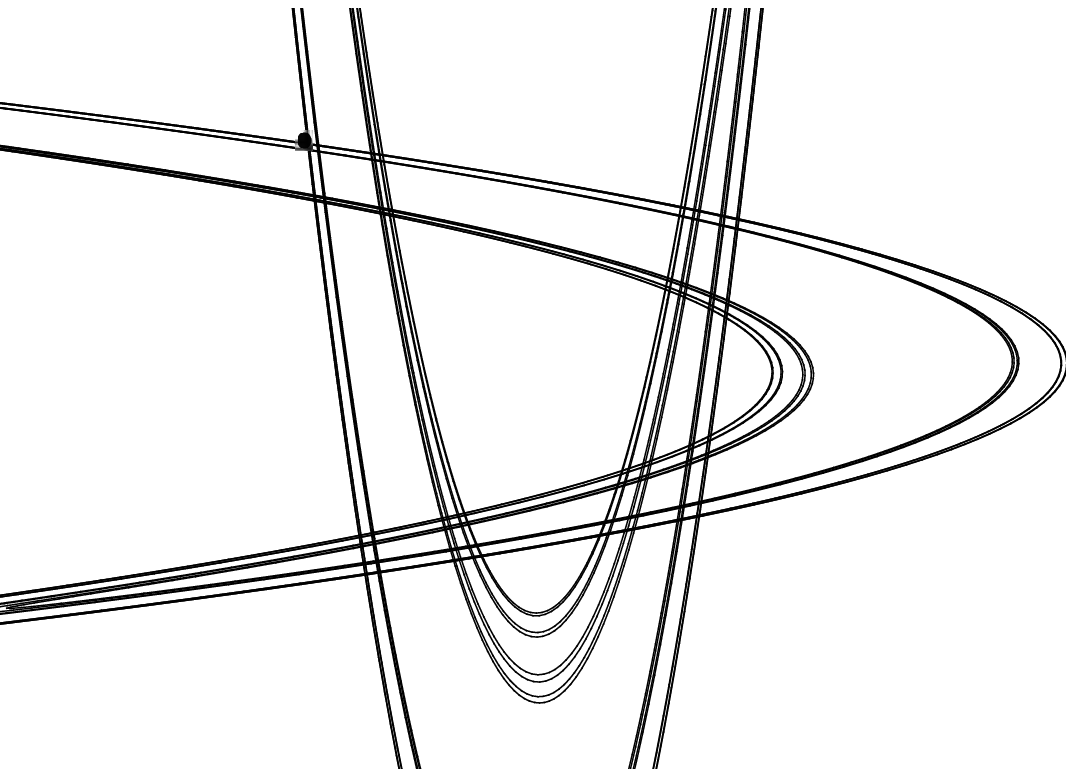}   }

\centerline{Stable/unstable manifolds for the Horseshoe Map}

\centerline{ $f_{c,a}(x,y)= (c + a y - x^2, -x)$ with $c=6.0$, $a = 0.8$}

\centerline{Figure 7.}
\medskip
\noindent We summarize some dynamical properties of the horseshoe map:

\item  {$\bullet$}   All points outside $\overline{W^s(x_0)} \cup \overline {W^u(x_0)}$ escape to infinity.
\item {$\bullet$}  $\overline{W^s(x_0)} \cap \overline {W^u(x_0)} \cong {\rm Cantor\ set}\times{\rm Cantor\ set}$
\item  {$\bullet$} Saddle (periodic) points are dense in ${\rm Cantor\ set}\times{\rm Cantor\ set}$
\item {$\bullet$}  Dynamics on ${\rm Cantor\  set}\times{\rm Cantor\ set}$ is conjugate to $\{0,1\}^{\Bbb Z}$
\item  {$\bullet$} The closure $\overline{W^s(z_0)}$ is the same for all saddles $z_0$.

\medskip
In the 1-dimensional case, the stable set is $W^s(x_0)=\bigcup_{n\ge0}f^{-n}(x_0)$.  This is an infinite set whose closure contains the Julia set.  If we start with an element $x_0$ in the Julia, then $f^{-n}(x_0)$ will be contained in the Julia set, and it will fill it out as $n\to\infty$.  By the Theorem at the beginning of this section, we may think of the invariant measure $\mu$ as describing how the sets $W^s(x_0)$ accumulate.  We want to do something similar in the 2-dimensional case.  Suppose that $\gamma$ is an oriented curve in ${\Bbb R}^2$, or if $\Gamma$ is a complex submanifold of ${\Bbb C}^2$.  If $\gamma$ has locally bounded length (or if $\Gamma$ has locally bounded are), we can consider the current of integration $[\gamma]$, which acts on test forms by $\varphi\mapsto \int_\gamma \varphi$.  In the case of the horseshoe, we cannot define the current of integration $[W^s(x_0)]$ directly, but we may construct a current by taking the average over arcs of stable manifolds in ${\cal W}^s$.  Thus, instead of an invariant measure, we have a family of transversal measures, which assign mass to families of stable arcs.   When we map ${\cal W}^s$ by $f^{-1}$, each individual stable manifolds is stretched, but the action on the set of transversal measures is very much like the Bernoulli trials that we saw in the 1-dimensional case.  Thus $f^{-1}$ maps one stable manifold to another, mixing them in the same chaotic way we saw in dimension 1.

\medskip
\epsfxsize=4.05in
\centerline{\epsfbox{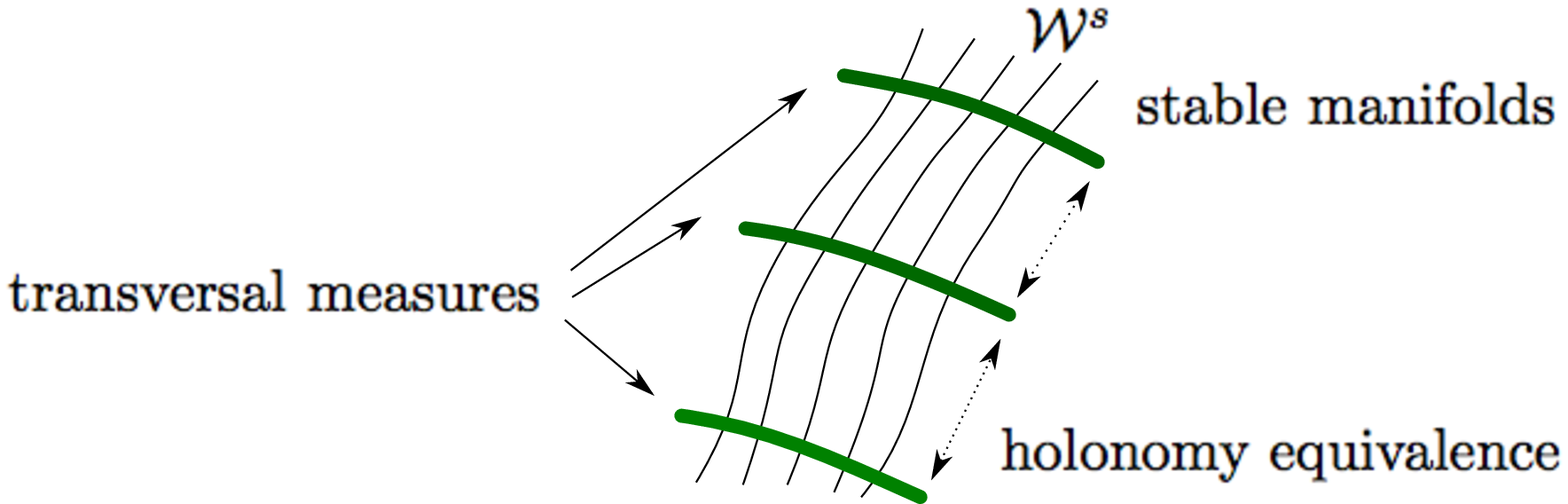}   }

\centerline{Mapping by $f^{-1}$ acts on the family of transversal measures}

\centerline{Figure 8.}
\medskip

While this phenomenon is rather special for real mappings, such as the horseshoe, it happens quite generally in the complex case.  If $f$ is an automorphism of a complex surface with $\delta(f)>1$, then there is an invariant current $T^s$ which satisfies $f^*T^s = \delta T^s$.  This current always exists, but it does not always have the same elegant geometric laminarity that we have seen in the case of the horseshoe.  The current $T^s$ can be sliced by a complex disk, and such slices $T^s|_D$ will serve as a family of ``transversal measures''.  If we work in a more measure-theoretic sense, we can assign a laminar structure to $T^s$, and such currents have been useful in studying the dynamics of $f$.  With Misha Lyubich and John Smillie, we have written a series of papers on the subject.  See  [D1--3] for a more recent treatment.

\bigskip\noindent {\bf \S6. What are the compact surfaces  $X$ which carry an $f\in Aut(X)$ with $\delta_1(f)>1$? }  This is answered in large part by a theorem of Cantat:
\proclaim Theorem [C1].  Suppose that $X$ is a compact complex manifold of (complex) dimension 2, and there is an automorphism $f$ of $X$ with $\delta(f)>1$. Then $X$ is a blowup of one of the following cases:
\item {$\bullet$}   $X$ is a torus, and $f$ is a ``standard'' torus automorphism.
\item {$\bullet$}   $X$ is a $K3$ surface, or a finite quotient of one of these.
\item {$\bullet$}   $X$ is a rational surface

The proof relies heavily on the Kodaira classification of surfaces.  It leads immediately to the question:  
{\it  What are the  surfaces and maps that actually occur for $K3$ or rational surfaces? }   The $K3$ surfaces or rational surfaces which can carry nontrivial automorphisms are quite special and not easy to find.

The set of all $K3$ surfaces has dimension 20, and the set of the ones that carry automorphisms has smaller dimension.  On the other hand, by [BK4], rational surfaces with automorphisms with $\delta>1$ can occur in families of arbitrarily large dimension.  This leads us to expect that a more interesting variety of dynamical behaviors will be found within the class of rational surfaces, so we concentrate on them, rather than on $K3$ surfaces.

\proclaim Theorem (Nagata).  In the case of a rational surface, we may suppose that $\pi:X\to{\Bbb P}^2$ is an (iterated) blowup of the projective plane.

Thus one possible approach is to look for birational maps $g:{\Bbb P}^2\dasharrow {\Bbb P}^2$ which can be lifted to automorphisms, as was done in \S2--4.  Another approach which has been productive is to require that the automorphism have an invariant curve.  McMullen gave a ``synthetic'' approach to this question, and Diller [Di] has shown all the possibilities for rational surface automorphisms with invariant curves, and which arise from quadratic maps of ${\Bbb P}^2$.
Let us remark (see [BK3]) that while many maps have been constructed, starting from invariant curves, not all rational surface automorphisms have invariant curves.  

If $X$ is obtained from ${\Bbb P}^2$ by a sequence of $N$ blowups, then $H^2(X;{\Bbb Z}) = {\Bbb Z}^{N+1}$, and the intersection product makes $H^2(X;{\Bbb Z})$ isometrically isomorphic to the Lorentz space ${\Bbb Z}^{1,N}$, with signature $(1,N)$.  A Theorem of Nagata (see [Dol1]) says that  if $f$ is an automorphism of $X$, then $f^*$ must belong to a Weyl group $W_N$.  Thus the question arises:  {\it  What are all the elements of $W_N$ which can be realized by rational surface automorphisms? }  That is, if $g\in W_N$, are there a surface $X$ and an $f\in Aut(X)$ such that the map induced by $f^*$ on ${\Bbb Z}^{1,N}$ coincides with  $g$?  If this is the case, then $\delta(f)$ is equal to the spectral radius of $g$.  Uehara [Ue] has shown that for any $g\in W_N$, there exists $f\in Aut(X)$  such that $\delta(f)$ is the spectral radius of $g$.  In other words, the set of all dynamical degrees coincides with the spectral radii of elements of $W_N$.  On the other hand, Diller [Di] and Uehara [Ue]  have shown that certain elements of $W_N$ cannot be realized by automorphisms with invariant curves, leaving open the question: {\it What happens for maps without invariant curves?}  The set of all dynamical degrees of surface automorphisms is well ordered and has other properties (see [BC]).  

\bigskip\noindent {\bf \S7.  Pseudo-automorphisms }  We were led to pseudo-automorphisms by the problem of  finding all periodicities within the family $f_{a,b}$  of birational maps of ${\Bbb P}^3$:
$$f_{a,b}(x_1,x_2,x_3) = \left( x_2,x_3,{ a_0 + a_1 x_1 + a_2 x_2 + a_3 x_3 \over b_0 + b_1 x_1 + b_2 x_2 + b_3 x_3} \right) $$
{\it What are the parameters $a_j$, $b_j$, $0\le j\le 3$, such that $f^p={\rm identity}$ for some $p$?}
This question seems to have originated when the first of the period 8 maps below was found many decades ago by Lyness [Ly].  The second one was found by Cs\"ornyei and Laczkovich [CL].   The following is proved in [BK5], where the word ``nontrivial'' is also explained.   
\proclaim Theorem.   The only nontrivial periods that can appear are 8 and 12.  The maps are 
$$(x, y, z) \to \left(y,z,{(1 + y+ z)/x}\right) {\rm\ \   or\ \   }  \left(y,z,{(-1 -y+ z)/ x}\right){\rm\ \  period\ 8 }$$
$$(x, y, z) \to\left(y,z,{({\eta\over 1-\eta} +\eta y+ z)/(\eta^2 x)}\right)  {\rm\ \    period\ 12 }$$
where $\eta$ is a primitive cube root of $-1$.

The principle behind the method of proof is to consider the function $(a,b)\mapsto\delta(f_{a,b})$ on parameter space and identify the subvariety of parameters $\{(a,b): \delta(f_{a,b})=1\}$, since any possible periodic map must lie here.  The generic map in this family has dynamical degree equal to $\delta_*>1$, and the scan of parameter space yielded a number of cases where $1<\delta(f_{a,b})<\delta_*$.  In the most interesting cases, we blow up 2 points, then 2 lines, and  then an orbit of curves (11 in one case and 19 in another) and arrive at a map without exceptional hypersurfaces.  The first of these cases is from [BK5, Theorem 1]:
$$F_a: \ \ \ (x,y,z)\mapsto (y,z, (a + \omega y + z)/x) \eqno(7.1)$$
where $a\in{\Bbb C}$, $a\ne0$, and $\omega$ is a primitive cube root of unity.  After blowups, these maps are ``regularized'' to the point that they are almost automorphisms, except that there are a finite number of indeterminate curves, and points of these indeterminate curves are blown up to other curves.  These maps are pseudo-automorphisms, and it is easy to deal with the map $f^*$ on cohomology, but the existence of indeterminacy makes it tricky to analyze the pointwise dynamics. 

In connection with the method described above, it would be interesting to know more generally for a family $f_a$ of rational maps:  {\it  Is $a\mapsto \delta(f_a)$ lower semi-continuous?  Is $\{a:\delta(f_a)\le t\}$ always a subvariety?}  The answer is ``yes'' for birational surface maps (see [X, Theorem 1.6]).  

The indeterminacy locus ${\cal I}(f)$ of any rational map $f:X\dasharrow Y$ has codimension at least 2.  Thus if $H$ is any hypersurface, we may define the image (proper transform) of $H$ as the closure of $f(H-{\cal I}(f))$.  We say that $H$ is {\it exceptional} if the codimension of the proper transform of $H$ is $\ge 2$. 
A birational  $f:X\dasharrow X$ is a {\it pseudo-automorphism} if neither $f$ nor $f^{-1}$ has an exceptional hypersurface.

{\it Pseudo-automorphisms behave very much in the spirit of automorphisms, and we expand our search to include this richer source of interesting maps.}   In dimension 2, all pseudo-automorphisms are in fact automorphisms.
What happens in  dimension $>2$?
Given that the blowups of ${\Bbb P}^2$ have yielded interesting automorphisms, it makes sense to ask: 
 {\it  Are there 3-folds $X$ which are obtained as blowups of ${\Bbb P}^3$ and which carry automorphisms $f$ with $\delta(f)>1$?}  Of course, we would expect such automorphisms to exist only in very special cases.  
\proclaim Theorem [T].   If $X$ is obtained from ${\Bbb P}^3$ by blowing up points and curves satisfying a certain condition, and if $f$ is an automorphism of $X$, then  $\delta_1(f)=\delta_2(f)$.

\proclaim Theorem [BaC].  If $X$ is obtained from ${\Bbb P}^k$ by blowing up points, then any automorphism $f$ of $X$ satisfies $\delta_\ell(f)=1$ for all $\ell$.

The Cremona involution on 
${\Bbb P}^3$ is the cubic map given by 
$$J (x) =
[1/x_0:1/x_1:1/x_2:1/x_3] = [x_1x_2x_3:x_0x_2x_3:x_0x_1x_3:x_0x_1x_2] $$
which acts as an involution on the coordinate tetrahedron $e_j\leftrightarrow \Sigma_j$, $j=0,1,2,3$.  We now see a new phenomenon:  any non-vertex point of an edge of the tetrahedron is blown up by $J$ to the skew edge.   In Figure 9, this means that any non-vertex point of $\alpha$ will be blown up to the whole edge $\beta$.

Let $\pi:X\to {\Bbb P}^3$ be the blowup of ${\Bbb P}^3$ at $e_0$,  The 3 edges of the tetrahedron passing through $e_0$ are now separated as in Figure~10.  The induced map $J_X:X\dasharrow X$  maps $E_0\leftrightarrow \Sigma_0$, so $\Sigma_0$ is no longer exceptional.  The restricted map $J_X|_{E_0}:E_0\to \Sigma_0$ ``looks like'' the 2D map $J$ mapping ${\Bbb P}^2$ to itself: the black triangle inside $E_0$ is exceptional, and the dotted black line is mapped to the bold black dot.  The  edges of the tetrahedron (in green) are still indeterminate.

\medskip
\epsfxsize=1.35in
\centerline{\epsfbox  {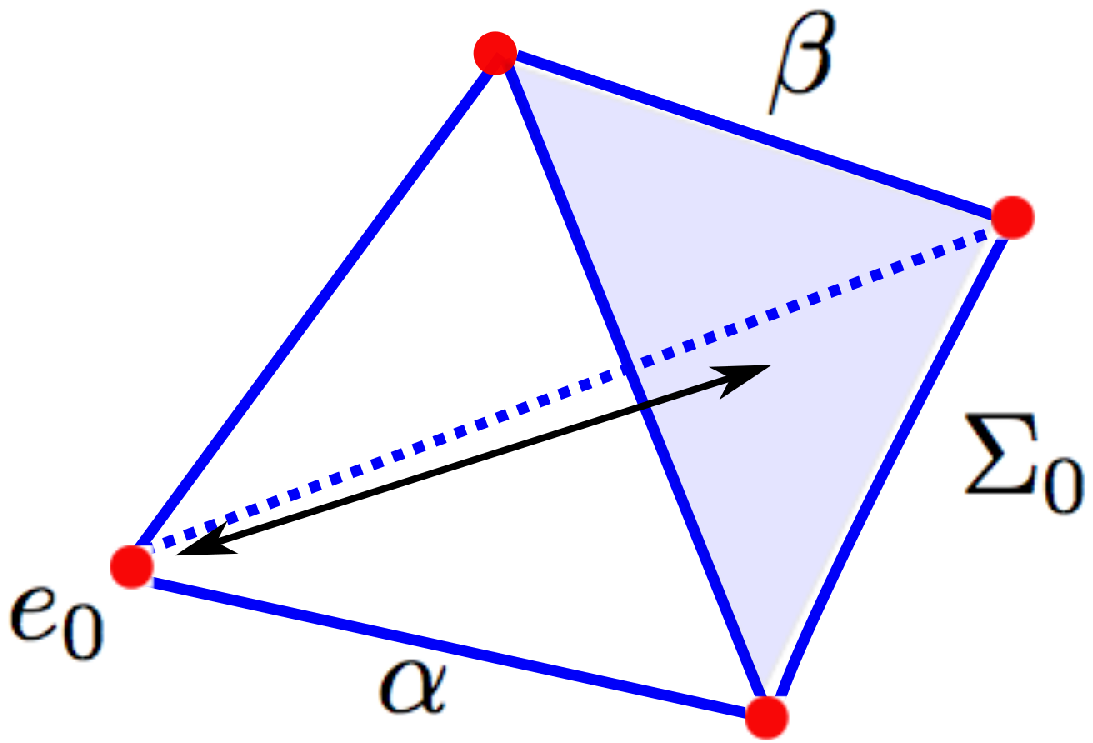}    }

\centerline{Cremona involution blows up any point of edge $\alpha$ to all of edge $\beta$.}

\centerline{Figure 9.}
\medskip

We may let $\pi:Y\to{\Bbb P}^3$ be the space obtained by blowing up all the vertices $e_j$, $j=0,1,2,3$, then the induced map $J_Y$ will be a pseudo-automorphism of $Y$.  In the space $Y$, the strict transforms of the edges of the tetrahedron are disjoint curves.  Let $\pi:Z\to Y$ be the space obtained by blowing up the strict transforms of the edges of the tetrahedron.  Thus $J_Z$ is an automorphism of $Z$, a very simple instance of the Theorem of de Fernex and Ein mentioned above.  

On the other hand, the example of $F_a$ in (7.1) illustrates why this procedure may encounter difficulties if the map is not periodic;  in [BK5, Theorem 4] it is shown that $F_a$ is not birationally conjugate to an automorphism.  Since blowing up or down is a birational operation, $F_a$ cannot be turned into an automorphism by any sort of blowing up procedure.  A heuristic explanation for this difficulty in the case of $F_a$ is that we will need to blow up an orbit of curves, but the curves in the orbit are not pairwise disjoint.

\medskip
\epsfxsize=2.05in
\centerline{\epsfbox  {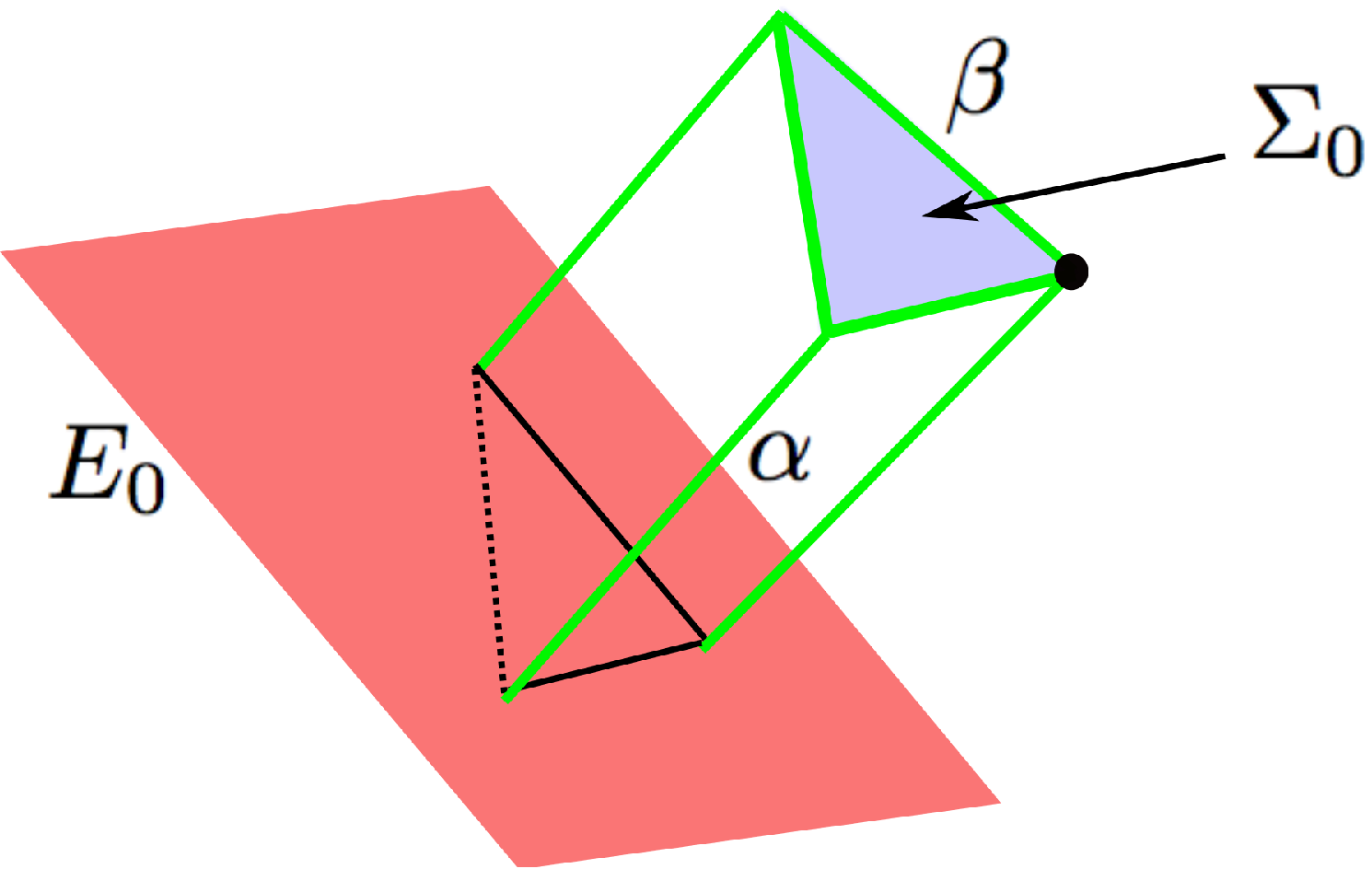}    }

\centerline{Cremona involution blows up any point of edge $\alpha$ to all of edge $\beta$.}

\centerline{Figure 10.}

\medskip\noindent{\bf \S8. Intermediate dynamical degrees}.  
Let $X$ be a manifold of dimension $k$, and let $f$ be an automorphism of $X$.  We measure the complexity of $f$ in terms of dynamical degrees.  We assume that $X$ is K\"ahler, so the action $f^*$ on cohomology respects the $(p,q)$ bi-gradation $H^*=\oplus H^{p,q}$.  If $1\le \ell\le k$, then the $\ell$-th dynamical degree is defined as the exponential rate of growth of the induced map $f^*$ on $H^{\ell,\ell}(X)$:
$$\delta_\ell(f):=\limsup_{n\to\infty}||f^{n*}:H^{\ell,\ell}(X)\to H^{\ell,\ell}(X)||^{1/n}$$
This is independent of the choice of norm on $H^{\ell,\ell}$. For an automorphism, we have $(f^n)^*=(f^*)^n$, so $\delta_\ell$  is the same as the spectral radius of  $f^*$ acting on $H^{\ell,\ell}(X)$.  It is also the same as the spectral radius of the restriction of $f^*$ to $H^{2\ell}(X)$.  As in \S2, an equivalent definition is the exponential rate of growth of the integrals $\int \beta^{k-\ell}\wedge (f^n)^*(\beta^{\ell})$ as $n\to\infty$.

In all cases, $\delta_\ell\ge1$, $\delta_0=1$, and $\delta_k={\rm deg}_{\rm top}$ is the topological or mapping degree of $f$.  Since $f$ is invertible, we have $\delta_k=1$.  By duality, $\delta_\ell(f) = \delta_{k-\ell}(f^{-1})$.  For the intermediate degrees, $\ell\mapsto \log(\delta_\ell)$ is concave.  Thus, if $\delta_1=1$, then $\delta_\ell=1$ for all $0\le \ell\le k$; and if $\delta_1>1$, then $\delta_\ell>1$ for all $1\le \ell\le k-1$. 

Let us note that if $f:X\dasharrow Y$ is merely rational, then there is still a well-defined linear map $f^*:H^*(Y)\to H^*(X)$.  
In dealing with rational maps, birational conjugacy is a natural sense of equivalence.
It was shown (see [DF] and [DS2]) that $\delta_\ell$ is an invariant of birational conjugacy.  However, the topological entropy is not a birational invariant (see [G]), but  [DS2] gives an inequality:
$${\rm entropy}_{\rm top}(f)\le \max(\log(\delta_1),\dots,\log(\delta_k))$$

The only general class of non-holomorphic, rational maps for which the intermediate degrees has been computed is the case of monomial maps, and the following is nontrivial:
\proclaim Theorem  [Lin], [FW].  Let $A=(a_{p,q})$ be an integer matrix of size $k\times k$, and let 
$$f_A(x) = x^A=(x_1^{a_{1,1}}\cdots x_k^{a_{1,k}},\dots,x_1^{a_{k,1}}\cdots x_k^{a_{k,k}})= \left(\prod_{q}x_q^{a_{1,q}},\dots,\prod_{q}x_q^{a_{k,q}} \right)$$ 
be the associated monomial map.  Then for each $\ell\ge 1$,  $\delta_\ell(f_A)$ is the spectral radius of the $\ell$-th exterior power ${\wedge^\ell}A$ of $A$.  Equivalently, if $|\lambda_1|\ge |\lambda_2|\ge\cdots\ge |\lambda_k|$ are the eigenvalues of $A$, then $\delta_\ell = |\lambda_1\lambda_2\cdots\lambda_\ell|$.

Intermediate degrees have also been determined by [A] and  [KR].  We have the problem:  
 {\it Determine $\delta_\ell(f)$ for $1<\ell<k$ for other nontrivial $f$.}

\bigskip\noindent{\bf \S9.  Existence of pseudo-automorphisms of blowups of ${\Bbb P}^k$ }  
We will try to find pseudo-automorphisms of the form $L\circ J$ for some $L$ which is nonsingular $(k+1)\times (k+1)$ matrix, and thus a linear automorphism  of ${\Bbb P}^k$.
The exceptional locus consists of the hyperplanes $\Sigma_j$, which are mapped:
$$f:=L\circ J \ : \ \Sigma_j\to L_j$$
where $L_j$ denotes the point of ${\Bbb P}^k$ defined by the $j$th column of the matrix $L$.  We will have a pseudo-automorphism if 
$$L_j\mapsto f(L_j)\mapsto\cdots\mapsto f^{m_j}(L_j)=e_{\sigma_j}, \ \ \ f^\ell(L_j)\notin \bigcup_i \Sigma_i$$
where $\sigma$ is a permutation of $\{0,\dots,k\}$.  Let $\pi:X\to{\Bbb P}^k$ be the blowup of the orbits of the $L_j$.  The induced map $f_X:=\pi^{-1}\circ f\circ \pi$ will be a pseudo-automorphism of $X$.  In [BK1, Theorem~A.1] a formula is given for the characteristic polynomial of $f_X$:
\proclaim Theorem. The polynomial  defining the dynamical degree $\delta_1(L\circ J)=\delta_1(m_1,\dots,m_k,\sigma)$ is given by an explicit formula involving the orbit lengths $m_j$ and the permutation $\sigma$.

Thus the lengths $m_j$ and the permutation $\sigma$ specify the dynamical degree that will be produced.  As a practical matter, however, the strategy given above for finding $L$ is not  feasible because it involves solving equations of very high degree in many variables.  The relevant computations are possible, however, if we assume that all the centers of blowup lie in an invariant curve. The existence of automorphisms of blowups of ${\Bbb P}^2$  with invariant curves was studied by McMullen,  Diller and Uehara.   Perroni and Zhang brought the method of McMullen from dimension 2 to higher dimension and gave the abstract existence of a map with an invariant curve.  
\proclaim Theorem [PZ].  For all $k\ge2$ and $d\ge1$ there exist infinitely many manifolds $X$ obtained by blowing up points on $({\bf P}^k)^d$ such that $X$ carries a pseudoautomorphism with $\delta>1$.

With Diller and Kim, we were motivated by the desire to see the Perroni-Zhang maps more concretely, and for this we used the method of Diller [Di].
Let us consider a parametrized curve $\psi: {\Bbb C}\to {\cal C}\subset {\Bbb P}^k$.  We say that ${\cal C}$ {satisfies a group law} if the following holds:  For each hyperplane $H\subset {\Bbb P}^k$, the set of all solutions (with multiplicity) $t_1,\dots,t_N$ of $\psi(t_i)\in H$ satisfies  $\sum t_i = 0$.  There are several cases of curves with group law; all the curves we work with have degree $k+1$.  For instance, there is ${\cal C}_1:=\psi({\Bbb C})$, which is the image of  $t\mapsto \psi(t)=[1:t:t^2:\cdots:t^{k-1}:t^{k+1}]$ which is irreducible and has a cusp.  There is also ${\cal C}_2:=\psi({\Bbb C})$, which is the image of the set-valued map $t\mapsto \psi(t) =\{[1:t:\cdots:t^k], [0:\cdots:0:(-1)^{k-1}:t]\}$ and is the  union of rational normal curve and a tangent line.  Both of these curves are singular at $[0:\cdots:0:1]=\psi(\infty)$.

If there is an invariant curve, then the points which will be blown up are of the form $\psi(t_j)$.  The problem of finding the centers of blowup is thus transformed to a problem of determining the points $\{t_i\}\subset{\Bbb C}$.  Using the group law on the curve, we  obtain:
\proclaim Theorem [BDK].  For most choices of orbit lengths $(m_0,\dots,m_k)$ and permutations $\sigma$, there is a matrix $L$ such that the space $\pi:X\to{\Bbb P}^k$ obtained by blowing up the orbits yields an induced pseudo-automorphism $f_X:X\dasharrow X$

This method applies also to products $({\Bbb P}^k)^d$.  For this we use the following variant of $J$.  We write a point of  $({\Bbb P}^k)^d$ as
$(x,y^{(1)},\dots,y^{(d-1)})$ and set
$$x/y^{(j)}:= (x_0/y^{(j)}_0,\dots,x_k/y^{(j)}_k), \ \ \ J(x,y^{(1)},\dots,y^{(d-1)}) = (1/x,x/y^{(1)},\dots,x/y^{(d-1)})$$
A linear map $L\in Aut(({\Bbb P}^k)^d$ has the form $L\circ\tau$ where $L:=(L_{1},\dots,L_{d})$,  $L_j\in Aut({\Bbb P}^k)$, and $\tau$ is a permutation of the factors ${\Bbb P}^k$.

This Cremona involution is discussed by Dolgachev [Dol2,3]  and Mukai [Muk] in connection with Weyl groups $W(p,q,r)$ which have $T$-shaped Coxeter diagrams.  In this case, $p=k+1$, $q=d+1$, and $r$ is the number of blowups $+1$ forming the space $\pi:
X\to ({\Bbb P}^k)^d$, and $W(p,q,r)\subset GL(H^2(X;{\Bbb Z}),{\Bbb Z})$.  It is shown in [BDK] that if a map of the form $f:=L\circ\tau\circ J$ is a pseudo-automorphism of a space $X$ obtained by blowing up  $({\Bbb P}^k)^d$, then the induced map $f^*$ belongs to the Weyl group $W(p,q,r)$.  More generally: {\it Let $X$ be a blowup of  $({\Bbb P}^k)^d$, and let ${\cal G}$ denote the group generated by $J$ and the  linear maps $Aut(({\Bbb P}^k)^d$.  If $f\in {\cal G}$ induces a pseudo-automorphism of $X$, does if follow that $f^*\in W(p,q,r)$?  } 

The simplest element of a Coxeter group with spectral radius $>1$ is the Coxeter element.  The results of [PZ] and [BDK] also apply to the Coxeter elements of $W(p,q,r)$.
For instance, let us consider the case with orbit lengths $(1,\dots, 1, n)$, and the permutation $\sigma=(0\,1\,2\, \dots k)$ is cyclic.  In this case the action of $f^*$ on $Pic(X)$ corresponds to the Coxeter element of the Weyl group $W(k+1,2, n+k+1)$.  The existence of such maps, representing the Coxeter  element, was given by Perroni and Zhang.  In [BDK] it is shown that we may write such a map as $f:=L\circ J$, where $L$ has the form:
$$L =\pmatrix{ 0& 0 & 0 & \dots & 0 & 1\cr
\beta_1 & 0 & 0 & \dots &0 & 1-\beta_1\cr
0 & \beta_2 & 0  & \dots & 0 & 1-\beta_2\cr
\vdots &&&&& \vdots \cr
0 & 0 &  0 &  \dots & \beta_k & 1-  \beta_k }$$
with $\beta_i$ being a rational function of $\delta$, any root of the minimal polynomial $\chi_n$ which gives the dynamical degree of $L\circ J$.  Different choices of invariant curve lead to different rational expressions for $\beta_i$.  

In dimension 2, it appears that the majority of automorphisms $f=L\circ J$ of blowups of ${\Bbb P}^2$ do not have invariant curves.  {\it Is it the case that most pseudo-automorphsms of $({\Bbb P}^k)^d$ do not have invariant curves?}   
\medskip
\noindent{\bf \S10.  Cohomological hyperbolicity. }
A map $f:X\dasharrow X$ is said to be {\it cohomologically hyperbolic} if there is a unique $1\le p\le k-1$ such that $\delta_p(f)$ is maximal.  In this case, the maximal growth occurs uniquely in bidegree $(p,p)$, which corresponds to codimension $p$.  As was noted in [BDK], we have $\delta_1=\delta_{k-1}$ for all maps $f=L\circ J$ which are pseudo-automorphisms of point blowups of ${\Bbb P}^k$.  Thus $f$ is not cohomologically hyperbolic when $k=3$. {\it  What are the intermediate dynamical degrees $\delta_\ell(L\circ J)$ when $k>3$?  Can $L\circ J$ be cohomologically hyperbolic?}

\proclaim Theorem  [DS1].  If $f$ is a cohomologically hyperbolic automorphism, then there are invariant currents $T^{s/u}$, and these may be used to form an invariant measure $\mu$ with interesting dynamical properties.  

  Guedj [G2] has conjectured that in the presence of cohomological hyperbolicity, the basic ergodic properties of 2-dimensional maps should carry over to higher dimension.    If $f$ is not cohomologically hyperbolic, then it is not clear to what extent a result like this would remain valid, and it is not clear what approach will reveal the dynamics of such maps.  It would be helpful if there could be an invariant fibration which would allow us to somehow study the dynamics with lower-dimensional objects and techniques.   In dimension~2, cohomological hyperbolicity fails exactly when $\delta(f)=1$.  In this case, Diller-Favre have shown:
\proclaim Theorem [DF].  If $f:X\dasharrow X$ is a bimeromorphic surface map with $\delta(f)=1$, then there is an invariant fibration.

Let $f:X\dasharrow X$ be a meromorphic map.  Suppose that there is a dominant, meromorphic map $\phi:X\dasharrow Y$ and a meromorphic map $g:Y\dasharrow Y$ such that $0<{\rm dim}(Y)<{\rm dim}(X)$, and $g\circ \phi = \phi\circ f$.  In other words, $\phi$ gives a meromorphic semiconjugacy from $(f,X)$ to $(g,Y)$.  In this case, the sets $\{\phi={\rm const}\}$ form an {\it invariant fibration}. 

In the presence of an invariant fibration, there is a {\it dynamical degree on the fiber}, written $\delta_j(f|\phi)$ (defined in [DN]), and it is related to the other dynamical degrees by:
\proclaim Theorem  [DN, DNT].  Suppose that the map $f$ has an invariant fibration as above.  Then $$\delta_p(f) = \max_{\max\{0,p-k+\ell
\}\le j\le \min\{p,\ell\} } \delta_j(g)\delta_{p-j}(f|\phi)$$

As a consequence of this, one can show:
\proclaim Theorem.
If $X$ is a 3-fold, and $f:X\dasharrow X$ is a birational map with an invariant fibration, then $\delta_1=\delta_2$.  In this case, $f$ is not cohomologically hyperbolic.

The possibilities for invariant fibrations in the automorphisms of tori are discussed in [OT1,2].  Guedj also conjectured: {\it  If $f$ is not cohomologically hyperbolic, then  $f$ has an invariant fibration, or at least an invariant foliation.}  In dimension 2, 
[KPR] have given a counter-example with $\delta_1=\delta_2= {\rm deg}_{\rm top}=2$, and thus is non-invertible.

If ${\rm dim}(X)=3$, and if $f:X\dasharrow X$ is birational, then the condition that $f$ is not cohomologically hyperbolic is equivalent to the condition that  $\delta_1(f)=\delta_2(f)$.  We now give a 3-dimensional invertible counterexample which is invertible.  Set
$$ L =\pmatrix{
0 & 0 & 0 & 1 \cr
1 & 0 & 0 &  a \cr
0 & 1 & 0  & 0 \cr
0 & 0 & 1 & c  
}$$
with $a,c\in{\Bbb C}$ such that 
$$n a^2 + (n+1)ac + n c^2=0$$ 
for some $n\ge2$, and let $J(x) = [1/x_0:\cdots:1/x_3]$ be the usual Cremona involution on ${\Bbb P}^3$.

\proclaim Theorem  [BCK].  For $n\ge 2$, we set $f_{a,c} := L\circ J$.  The dynamical degrees are $\delta_1(f)=\delta_2(f)>1$.  There is no (singular) foliation of dimension 1 or 2 which is invariant under $f_{a,c}$.   In particular, there is no invariant (singular) fibration.

{Structure of example}
$$f(\Sigma_0)=e_1:=[0:1:0:0], \ \  f(\Sigma_1)=e_2, \ \  f(\Sigma_2)=e_3$$
$$f(\Sigma_3)=p:=[1:a:0:c]$$

\proclaim Theorem. 
Let $Y$ denote ${\Bbb P}^3$ blown up at the points $e_1$ and $e_3$.  Then the induced map $f_Y$ is a dominant map of an invariant 4-cycle of surfaces:  $$\Sigma_0\to E_1\to \Sigma_2\to E_3\to \Sigma_0$$  
The orbit of the exceptional image point $p$ is inside this invariant set.
\medskip
\epsfxsize=2.05in
\centerline{\epsfbox  {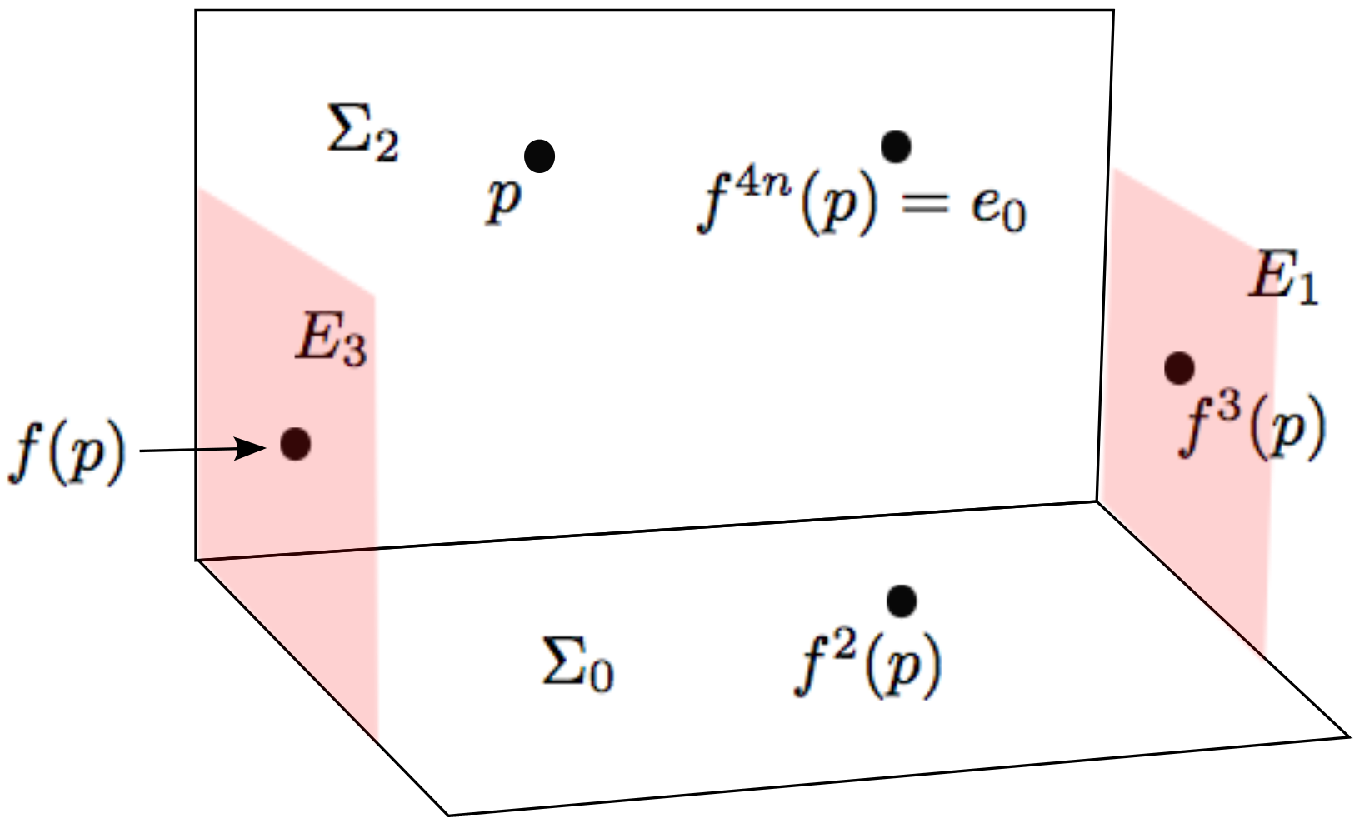}    }

\centerline{Construction of $Y$.}

\centerline{Figure 11.}

\proclaim Theorem.
If $\ n a^2 + (n+1)ac + n c^2=0$, then the $f_Y$ orbit of $q$ lands on the point $e_0$.  Let $X$ denote the space obtained by blowing up the $4n+2$ points $p$, $f_Y(p)$, \dots, $f_Y^{4n}p = e_0$, and $e_3$.  Then the induced map $f_X$ is a pseudo-automorphism. 

One difference between these maps and the [BDK] maps is that that there are two ``levels'' of blowup, and in fact none of the [BDK] maps is birationally conjugate to any of the [BCK] maps.  The invariant 4-cycle of surfaces $\Gamma:=\Sigma_0\cup E_1\cup\Sigma_2\cup E_3$ plays an important role in understanding $f_{a,c}$.  \proclaim Theorem.  For $g:=f^4|_{\Sigma_0}$, the dynamical degree satisfies $\delta_1(g)>1$, but it is not a Salem number.  Thus $g$ is not birationally conjugate to a surface automorphism.  

We note that the maps $F_a$ of (7.1) were analyzed by means of an invariant 8-cycle of surfaces with similar properties.  We may use [BK5, Theorem 1.5] to conclude that the maps  in [BCK] are not birationally conjugate to an automorphism of a 3-dimensional manifold.  

%The divisor $\Gamma$ represents the anti-canonical class, and it spans the 1-eigenspace of $f^*_X$ acting on $H^{1,1}(X)$.  Every $f_X$-invariant divisor is a multiple of $\Gamma$.  By studying $g^*$, we may show:
%
%\proclaim Theorem.  There is a unique $g$-invariant curve on $\Gamma$.
%
%{One idea in the proof of the BCK Theorem}
%Suppose that ${\cal F}$ is an invariant foliation of codimension 1.  Then the singular locus of ${\cal F}$ is an algebraic variety of dimension at most 1.
%If ${\cal F}$ ``crosses'' one of the indeterminate curves, then the image of that curve under $J$ must be in the singular set of $J$.
%
%
%
%\medskip
%\epsfxsize=4.05in
%\centerline{\epsfbox  {figures/singularity.eps}    }
%
%\centerline{Figure 12.}
%
%
%
%
%The image of one of these singular curves must lie inside the invariant 4-cycle of surfaces.  It is trapped inside this 4-cycle and has infinite orbit there.  This contradicts the fact that it must lie inside the 1-dimensional singular locus of ${\cal F}$.
%

\bigskip
\centerline{\bf References}

\item{[A]}   E. Amerik,  A computation of invariants of a rational self-map. Ann. Fac. Sci. Toulouse Math. (6) 18 (2009), no. 3, 445--457.

\item{[BaC]}   Bayraktar, Turgay; Cantat, Serge Constraints on automorphism groups of higher dimensional manifolds. J. Math. Anal. Appl. 405 (2013), no. 1, 209--213. 

\item{[B]}  E. Bedford,  The dynamical degrees of a mapping. Proceedings of the Workshop Future Directions in Difference Equations, 3--13, Colecc. Congr., 69, Univ. Vigo, Serv. Publ., Vigo, 2011. 

\item{[BCK]}  Eric Bedford, Serge Cantat and Kyounghee Kim,   Pseudo-automorphisms with no invariant foliation, {\sl Journal of Modern Dynamics}, to appear.  {  arXiv:1309.3695} 

\item{[BD]}  Eric Bedford and Jeffrey Diller,  Real and complex dynamics of a family of birational maps of the plane: the golden mean subshift. Amer.\ J. Math.\ 127 (2005), no.\ 3, 595--646. 

\item{[BDK]}  Eric Bedford, Jeffrey Diller and Kyounghee Kim, Pseudoautomorphisms with invariant curves,  { arXiv:1401.2386} 

\item{[BK1]}  Eric Bedford and Kyounghee Kim,  On the degree growth of birational mappings in higher dimension. J. Geom. Anal. 14 (2004), no. 4, 567--596. 

\item{[BK2]}    Eric Bedford and Kyounghee Kim,  Periodicities in linear fractional recurrences: degree growth of birational surface maps. Michigan Math. J. 54 (2006), no. 3, 647--670. 

\item{[BK3]}  Eric Bedford and Kyounghee Kim,   Dynamics of rational surface automorphisms: linear fractional recurrences. J. Geom. Anal. 19 (2009), no. 3, 553--583. 

\item{[BK4]}  Eric Bedford and Kyounghee Kim,  Continuous families of rational surface automorphisms with positive entropy. Math. Ann. 348 (2010), no. 3, 667--688.

\item{[BK5]} E. Bedford and K. Kim,  Dynamics of (pseudo) automorphisms of 3-space: periodicity versus positive entropy. Publ. Mat. 58 (2014), no. 1, 65--119. 

\item{[BC]}   J\'er\'emy Blanc and Serge Cantat, 
    Dynamical degrees of birational transformations of projective surfaces.
    { arXiv:1307.0361}

\item{[Br]}  Hans Brolin, Invariant sets under iteration of rational functions.
Ark. Mat. 6 1965 103--144 (1965). 

\item{[C1]}  Serge Cantat,   Dynamique des automorphismes des surfaces projectives complexes.  C. R. Acad. Sci. Paris S\'er. I Math. 328 (1999), no. 10, 901--906.

\item{[C2]}  Serge Cantat, Quelques aspects des syst\`emes dynamiques polynomiaux: existence, exemples, rigidit\'e.  13--95, Panor. Synth\`eses, 30, Soc. Math. France, Paris, 2010. 

\item{[CL]}   M. Cs\"ornyei and M. Laczkovich, Some periodic and non-periodic recursions. Monatsh. Math. 132 (2001), no. 3, 215--236.

\item{[dFE]}  T. de Fernex and L. Ein, Resolution of indeterminacy of pairs. Algebraic geometry, 165--177, de Gruyter, Berlin, 2002.

\item{[Di]}   Jeffrey Diller,  Cremona transformations, surface automorphisms, and plane cubics. With an appendix by Igor Dolgachev. Michigan Math. J. 60 (2011), no. 2, 409--440.

\item{[DS1]}    T.-C. Dinh and N. Sibony, Green currents for holomorphic automorphisms of compact K\"ahler manifolds. J. Amer. Math. Soc. 18 (2005), no. 2, 291--312.

\item{[DS2]}   T.-C. Dinh and N. Sibony, Une borne sup\'erieure pour l'entropie topologique d'une application rationnelle. Ann. of Math. (2) 161 (2005), no. 3, 1637--1644. 

\item{[Do1]}  Igor Dolgachev,   Reflection groups in algebraic geometry. Bull. Amer. Math. Soc. (N.S.) 45 (2008), no. 1, 1--60. 

\item{[Do2]}  Igor Dolgachev,   Infinite Coxeter groups and automorphisms of algebraic surfaces. The Lefschetz centennial conference, Part I (Mexico City, 1984), 91--106, Contemp. Math., 58, Amer. Math. Soc., Providence, RI, 1986. 

\item{[Do3]}  Igor Dolgachev,   Weyl groups and Cremona transformations. Singularities, Part 1 (Arcata, Calif., 1981), 283--294, Proc. Sympos. Pure Math., 40, Amer. Math. Soc., Providence, RI, 1983.

\item{[D1]}  R. Dujardin, Structure properties of laminar currents on ${\bf  P}\sp 2$. J. Geom. Anal. 15 (2005), no. 1, 25--47.

\item{[D2]} R. Dujardin, Sur l'intersection des courants laminaires.   Publ. Mat. 48 (2004), no. 1, 107--125.

\item{[D3]} R. Dujardin, Laminar currents in ${\bf P}\sp 2$. Math. Ann. 325 (2003), no. 4, 745--765.

\item{[FW]} C. Favre and E. Wulcan,  Degree growth of monomial maps and McMullen's polytope algebra,  Indiana Univ. Math. J. 61 (2012), no. 2, 493--524.  { arXiv:1011.2854}

\item{[FS1]}  J.E. Forn\ae ss and N. Sibony,   Complex dynamics in higher dimensions. Notes partially written by Estela A. Gavosto. NATO Adv. Sci. Inst. Ser. C Math. Phys. Sci., 439, Complex potential theory (Montreal, PQ, 1993), 131--186, Kluwer Acad. Publ., Dordrecht, 1994.

\item{[FS2]}  J.E. Forn\ae ss and N. Sibony,  Complex dynamics in higher dimension. II. Modern methods in complex analysis (Princeton, NJ, 1992), 135--182, Ann. of Math. Stud., 137, Princeton Univ. Press, Princeton, NJ, 1995.

\item{[FLM]}  A. Freire, A. Lopes and R. Ma\~n\'e, 
An invariant measure for rational maps.
Bol. Soc. Brasil. Mat. 14 (1983), no. 1, 45--62. 

\item{[G1]}  V. Guedj,   Entropie topologique des applications m\'eromorphes.  Ergodic Theory Dynam.\ Systems 25 (2005), no. 6, 1847--1855.

\item{[G2]}  V. Guedj,   Propri\'et\'es ergodiques des applications rationelles,  {\sl Quelques aspects des syst\`emes dynamiques polyn\^omiaux}, 97--202, Panor.\ Synth\`eses, 30, Soc.\ Math.\ France, Paris, 2010.
{ arXiv:math/0611302}

\item{[KPR]}    Scott R. Kaschner, Rodrigo A. P\'erez, Roland K.W. Roeder, Examples of rational maps of ${\Bbb CP}^2$ with equal dynamical degrees and no invariant foliation, {  arXiv:1309.4364} 

\item{[KR]}   Sara Koch and Roland Roeder,  Computing dynamical degrees, {  arXiv:1403.5840} 

\item{[Lin]}  Jan-Li Lin,  Pulling back cohomology classes and dynamical degrees of monomial maps. Bull.\ Soc.\ Math.\ France 140 (2012), no. 4, 533--549 (2013). 

\item{[Ly]}  R.C. Lyness, Notes 1581,1847, and 2952, Math.\ Gazette {\bf 26} (1942), 62, {\bf 29} (1945), 231, and {\bf 45} (1961), 201.

\item{[L]}  M. Lyubich,  Entropy properties of rational endomorphisms of the Riemann sphere.
Ergodic Theory Dynam. Systems 3 (1983), no. 3, 351--385.

\item{[M]}  Curtis T. McMullen,  Dynamics on blowups of the projective plane. Publ. Math. Inst. Hautes ƒtudes Sci. No. 105 (2007), 49--89.

\item{[Muk]}  Shigeru Mukai,  Geometric realization of $T$-shaped root systems and  counterexamples to Hilbert's tenth problem.  In {\sl Algebraic transformation groups and algebraic varieties}, volume 132 of {\sl Encyclopedia  Math.\ Sci.}, pages 123--129.  Springer, Berlin.  2004

\item{[OT1]}  Keiji Oguiso and Tuyen Trung Truong, Salem numbers in dynamics of K\"ahler threefolds and complex tori,  {  arXiv:1309.4851}

\item{[OT2]}  Keiji Oguiso and Tuyen Trung Truong,    Explicit Examples of rational and Calabi-Yau threefolds with primitive automorphisms of positive entropy, {  arXiv:1306.1590  }

\item{[PZ]}  Fabio Perroni and De-Qi Zhang, Pseudo-automorphisms of positive entropy on the blowups of products of projective spaces.\ Math.\ Ann.\ 359 (2014), no.\ 1-2, 189--209.

\item{[T]}  Tuyen Trung Truong,    On automorphisms of blowups of ${\bf P}^3$,  { arXiv:1202.4224}

\item{[Ueh]} T.  Uehara,  Rational surface automorphisms with positive entropy, 
{  arXiv:1009.2143}

\item{[X]}   Junyi Xie, Periodic points of birational maps on projective surfaces, {  arXiv:1106.1825} 

\bigskip\bigskip
\rightline{Stony Brook University}

\rightline{Stony Brook, NY 11794}

\rightline{ebedford@math.sunysb.edu}

\bye